\documentclass[11pt,reqno]{amsart}
\usepackage{amsmath}
\usepackage{txfonts}
\usepackage{amssymb}
 \usepackage{xcolor}
\usepackage{amsmath,amsthm}
\usepackage{CJK,amsfonts}
\usepackage[body={7.5in,8.75in}, top=1.2in, right=0.9in,left=0.9in]{geometry}
\usepackage{hyperref}%,%,backref,
\allowdisplaybreaks[4] \numberwithin{equation}{section}

\newtheorem{thm}{Theorem}

\newtheorem{cor}[thm]{Corollary}
\newtheorem{lem}[thm]{Lemma}

\def\squarebox#1{\hbox to #1{\hfill\vbox to #1{\vfill}}}

\newcommand{\be}{\begin{equation}}
\newcommand{\ee}{\end{equation}}
\newcommand{\bea}{\begin{eqnarray}}
\newcommand{\eea}{\end{eqnarray}}
\newcommand{\bd}{\begin{displaymath}}
\newcommand{\ed}{\end{displaymath}}
\begin{document}
\begin{CJK*}{GBK}{song}
\title[Connection formulae for trivariate $q$-polynomials]{Connection formulae for trivariate $q$-polynomials}
\author{    Sama Arjika${}^{1}$ and Zouha\"{\i}r Mouayn$^{\ast }$}
\dedicatory{\textsc}
\thanks{Address for correspondence: ${}^{1}$Department of Mathematics and Informatics, University of Agadez, Niger. ${}^2$Department of Mathematics, Faculty of Sciences
and Technics (M'Ghila)  BP. 523, B\'{e}ni Mellal, Morocco }
\thanks{Email:  rjksama2008@gmail.com}

\keywords{binomial coefficients; Cauchy polynomials; Connection formulae.  }

\thanks{2010 \textit{Mathematics Subject Classification}. 05A15,
11B65, 26A33, 33D15, 33D45, 33D60, 39A13.}

\begin{abstract}
In this short paper, we   establish connection formulae for trivariate $q$-polynomials.  
\end{abstract}

\maketitle

\section{Introduction, notations and definitions}
 The $q$-analysis  can be traced back to the earlier works of L. J. Rogers  \cite{RG94}. It has wideranging applications in the analytic number theory  and $q$-deformation of well-known functions \cite{GEA} as well as in the  study of solvable models in statistical mechanics \cite{RJB}. During the 80's the interest on this  analysis increased with quantum groups theory  with which models of $q$-deformed oscillators have been developed \cite{MJ}. The $q$-analogs of boson operators have been defined in \cite{ARS} where the corresponding wavefunctions were constructed  in terms of the continuous $q$-Hermite polynomials of Rogers and other polynomials. Actually, known models of $q$-oscillators are closely related with $q$-orthogonal polynomials.
 
 In this paper, we adopt the common   notation and terminology for basic hypergeometric  series as in Refs.  \cite{GasparRahman,Koekock}.    Throughout this paper, we assume that  $q$ is a fixed nonzero real or complex number and $|q|< 1$.  The $q$-shifted factorial  and its compact  factorial  are defined \cite{GasparRahman,Koekock}, respectively by:
\begin{equation}
(a;q)_0:=1,\quad  (a;q)_{n} =\prod_{k=0}^{n-1} (1-aq^k),   \quad (a;q)_{\infty}:=\prod_{k=0}^{\infty}(1-aq^{k})
\end{equation}
and 
 $ (a_1,a_2, \ldots, a_r;q)_m=(a_1;q)_m (a_2;q)_m\cdots(a_r;q)_m,\; m\in\{0, 1, 2\cdots\}$. 

 Here, in our present investigation, we are mainly concerned  
with the 
 Jackson-Hahn-Cigler (JHC)  $q$-addition 
 %\begin{center*}
\begin{equation}
(x\oplus _{q}y)^{n}:=\sum_{k=0}^{n}\begin{bmatrix} n\\ k \end{bmatrix}_q q^{\binom{k}{2}
}x^{n-k}y^{k}=x^n\left(\frac{-y}{x};q\right)_n\equiv P_{n}(x,-y),\: n=0,1,2,\cdots,\label{2.7}
\end{equation}
where the
%The JHC $q$-substruction is defined by 
%\begin{equation*}
%(x\ominus _{q}y)^n:=P_{n}(x,-y),\: n=0,1,\cdots,  
%\end{equation*}
Cauchy polynomials $P_n(x,y)$ as given
below (see \cite{Chen2003} and \cite{GasparRahman}):
\bea
\label{def}
 P_n(x,y):=(x-y)( x- qy)\cdots ( x-q^{n-1}y) =(y/x;q)_n\,x^n:=(x\ominus _{q}y)^n
\eea
 which has the following  generating function \cite{Chen2003}
\be
\label{gener}
\sum_{n=0}^{\infty} P_n(x,y)
\frac{t^n }{(q;q)_n} = 
\frac{(yt;q)_\infty}{(xt;q)_\infty}.
\ee
 The generating   function (\ref{gener}) is also the   homogeneous version  of the Cauchy identity or the $q$-binomial theorem   given by \cite{GasparRahman}
\be
\label{putt}
\sum_{k=0}^{\infty} 
\frac{(a;q)_k }{(q;q)_k}z^{k}={}_{1}\Phi_0\left[\begin{array}{c}a;
 \\\\
-;
 \end{array} 
q;z\right]= 
\frac{(az;q)_\infty}{(z;q)_\infty},\quad |z|<1,  
\ee
where the    basic or $q$-hypergeometric function    in the variable $z$ (see  Slater \cite[Chap. 3]{SLATER},  Srivastava and Karlsson   \cite[p. 347, Eq. (272)]{SrivastaKarlsson}   for details) is defined as:
 $$
{}_{r}\Phi_s\left[\begin{array}{r}a_1, a_2,\ldots, a_r;
 \\\\
b_1,b_2,\ldots,b_s;
 \end{array} 
q;z\right]
 =\sum_{n=0}^\infty\Big[(-1)^n q^{({}^n_2)}\Big]^{1+s-r}\,\frac{(a_1, a_2,\ldots, a_r;q)_n}{(b_1,b_2,\ldots,b_s;q)_n}\frac{ z^n}{(q;q)_n},
$$
 when $r>s+1$. Note that, for $r=s+1$, we have:
$$
{}_{r+1}\Phi_r\left[\begin{array}{r}a_1, a_2,\ldots, a_{r+1};
 \\\\
b_1,b_2,\ldots,b_r;
 \end{array} q; z\right]
 =\sum_{n=0}^\infty \frac{(a_1, a_2,\ldots, a_{r+1};q)_n}{(b_1,b_2,\ldots,b_r;q)_n}\frac{ z^n}{(q;q)_n}.
$$
Putting    $a=0$, the relation (\ref{putt}) becomes  Euler's identity  \cite{GasparRahman}
\be
\label{q-expo-alpha}
  \sum_{k=0}^{\infty} \frac{ z^{k}}{(q;q)_k}=\frac{1}{(z;q)_\infty},\quad |z|<1
\ee
and its inverse relation  \cite{GasparRahman}
\be
\label{q-Expo-alpha}
 \sum_{k=0}^{\infty}  \frac{(-1)^kq^{ ({}^k_2)
}\,z^{k}}{(q;q)_k}=(z;q)_\infty.
\ee
The relations (\ref{q-expo-alpha}) and (\ref{q-Expo-alpha}) satisfy
\begin{equation}
e_{q}(a)E_{q}(b)=e_{q}(a\oplus _{q}b),\quad e_{q}(a)E_{q}(-a)=1. \label{2.12}
\end{equation}
 
In  \cite{Mohammed}, Mohammed introduced the trivariate $q$-polynomials as
\begin{equation}
F_n(x,y,z;q)=(-1)^nq^{-({}^n_2)} \sum_{k=0}^n {\,n\,\atopwithdelims []\,k\,}_{q}(-1)^kq^{({}^k_2)}  P_{n-k}(y,x)z^k \label{3.1}
\end{equation}%
with the following generating function \cite[Theorem 2.6]{Mohammed}
\begin{equation}
\sum_{n=0}^{\infty } F_n(x,y,z;q)\frac{ (-1)^nq^{ ({}^n_2)}\,t^n}{(q;q)_n}=\frac{(xt,zt,q)_\infty}{(yt; q)_\infty}.  \label{3.6}
\end{equation}

In this paper, we shall  establish a  connection 
 formulae for trivariate $q$-polynomials.   The  connection formulae for  orthogonal polynomials are useful in mathematical analysis and also have applications in quantum mechanics, such as finding the relationships between the wavefunctions of some  potential functions used for describing physical and chemical properties in atomic and molecular systems \cite{SLD}.

In the rest, we   we establish some connection formulae and discuss some particular cases. 

\section{Main results}\label{section2}
In this section, we give the following  fundamental  theorem. 
\begin{thm}\label{thm1}
Let $f(x,y,z)$ be a three-variable analytic function in a neighborhood of $(x,y,z)=(0,0,0)\in\mathbb{C}^3.$ Then, $f(x,y,z)$ can be expanded in terms of  $F_n(x,y,z;q)$  if and only if $f$ satisfies the following $q$-difference equation:
\begin{align}
\label{asEZA}
 (q^{-1}x- y)\Big[f(x,y, z)-f(x,y,qz)\Big]
 =  z  [f(q^{-1}x,y, qz)-zf(x, qy, qz).
\end{align}
%then we have:
%\bea
%\label{shha}
%f(x,y, z)=\widetilde{L}(a,z;\theta_{xy})\Big\{f(a ,x,y, 0) \Big\}.
%\eea
\end{thm}
To determine if a given function is an analytic function in several complex variables, we often use the following Hartogs's Theorem. For more information, please refer to
Taylor \cite[p. 28]{Taylor} and Liu \cite[Theorem 1.8]{Liu}.
\begin{lem} 
  $\big[$Hartogs's Theorem \cite[p.15]{Gunning}$\big]$
\label{lemma1}
If a complex-valued function is holomorphic (analytic) in each variable separately in an open domain $D\subset\mathbb{C}^n$, then it is holomorphic (analytic) in $D.$
\end{lem}
\begin{lem}\cite[p. 5 Proposition 1]{Malgrange}
\label{lemma2}
If $f(x_1, x_2, . . . , x_k)$ is analytic at the origin $
(0, 0, . . . , 0)\in\mathbb{C}^k$, then, $f$ can be expanded in an absolutely convergent power series
\be
f(x_1, x_2, . . . , x_k)=\sum_{n_1,n_2,\cdots,n_k=0}^\infty \alpha_{n_1,n_2,\cdots,n_k} x_1^{n_1} x_2^{n_2}\cdots  x_k^{n_k}.
\ee

\end{lem}

\begin{proof}[Proof of Theorem \ref{thm1}] From the Hartogs's Theorem and the theory of several complex variables (see Lemmas \ref{lemma1} and \ref{lemma2}), we assume that
\be
\label{a120}
f(x,y, z)=\sum_{n=0}^\infty A_n(x,y)z^n.
\ee
Then,  substituting the above equation into (\ref{asEZA}), we have: 
\begin{align}
(q^{-1}x- y)\sum_{n=0}^\infty (1-q^n)A_n(ax,y)z^n 
=  \sum_{n=0}^\infty q^n [A_n(q^{-1}x,y)-A_n(x,qy) 
 ]z^{n+1}
\end{align}
Comparing coefficients of $z^n,\,n\geq 1$, we readily find that
\be 
(q^{-1}x- y)  (1-q^n)A_n(x,y) 
=   q^{n-1}(1-aq^{n-1})[A_{n-1}(q^{-1}x,y)-A_{n-1}(x,qy) 
].
\ee
After simplification, we get 
\bea
\label{sZA}
 A_n(x,y)=     \frac{q^{n-1} }{1-q^n} \theta_{xy}\Big\{ A_{n-1}(x,y)\Big\}.
\eea
By iteration, we gain
\bea
\label{ZddA}
  A_n(x,y)=   \frac{q^{({}^n_2)} }{(q;q)_n}\theta_{xy}^n\Big\{ A_0(x,y)\Big\}.
\eea
Letting $\displaystyle f(x,y,0)=A_0(x,y)=\sum_{n=0}^\infty \mu_np_n(y,x),$ we have
\bea
\label{dA}
   A_k(a,x)=    \frac{  q^{({}^k_2)} }{(q;q)_k}\sum_{n=0}^\infty \mu_n\frac{(-1)^k(q;q)_n }{(q;q)_{n-k}}p_{n-k}(y,x).
\eea
Replacing  (\ref{dA}) in (\ref{a120}), we have:
\bea
f(a,x,y)&=&\sum_{k=0}^\infty \frac{ (-1)^kq^{({}^k_2)} }{(q;q)_k}\sum_{n=0}^\infty \mu_n\frac{(q;q)_n }{(q;q)_{n-k}}p_{n-k}(y,x)z^k\cr
&=&\sum_{n=0}^\infty \mu_n\sum_{k=0}^n
\begin{bmatrix} n\\k \end{bmatrix}_{q}
  (-1)^kq^{({}^k_2)} p_{n-k}(y,x)z^k. 
\eea
On the other hand, if  $f(x,y,z)$ can be expanded  in term of $ F_n(x,y,z)$, we can verify that  $f(x,y,z)$ satisfies (\ref{asEZA}). The proof of the assertion (\ref{asEZA}) of Theorem \ref{thm1} is now completed. 
\end{proof}

\begin{thm}
\label{aprodpos}
Let $f(x,y,z)$ be a three-variable analytic function in a neighborhood of $(x,y,z)=(0,0,0)\in\mathbb{C}^3$. 
 
    If $f(x,y,z)$  satisfies the $q$-difference equation
\begin{align}
\label{sEZA}
 (q^{-1}x- y)\Big[f(x,y, z)-f(x,y,qz)\Big]
 =  z \Big[f(q^{-1}x,y, qz)-f(x, qy, qz)\Big]
\end{align}
then we have:
\bea
\label{shha}
f( x,y, z)= L(z\theta_{xy})\Big\{f(x,y, 0) \Big\}.
\eea

\end{thm}
\begin{proof}[Proof of Theorem \ref{aprodpos}]  From the theory of several complex variables \cite{Malgrange},  we   begin to solve the $q$-difference equation (\ref{sEZA}). First we may assume that
\be
\label{120}
f(x,y, z)=\sum_{n=0}^\infty A_n(x,y)z^n.
\ee
Then  substituting the above equation into (\ref{sEZA}), we have: 
\begin{align}
(q^{-1}x- y)\sum_{n=0}^\infty (1-q^n)B_n(ax,y)z^n 
=  \sum_{n=0}^\infty q^n [A_n(q^{-1}x,y)-A_n(x,qy) 
 ]z^{n+1}
\end{align}
Comparing coefficients of $z^n,\,n\geq 1$, we readily find that
\be 
(q^{-1}x- y)  (1-q^n)B_n(a,x,y) 
=   q^{n-1}(1-aq^{n-1})[B_{n-1}(a, q^{-1}x,y)-B_{n-1}(a, x,qy) 
].
\ee
After simplification, we get 
\bea
\label{sZA}
 A_n(x,y)=     \frac{q^{n-1} }{1-q^n} \theta_{xy}\Big\{ A_{n-1}(x,y)\Big\}.
\eea
By iteration, we gain
\bea
\label{ZddA}
  A_n(x,y)=   \frac{q^{({}^n_2)} }{(q;q)_n}\theta_{xy}^n\Big\{ A_0(x,y)\Big\}.
\eea
Now we return to calculate $A_0(x,y)$. Just taking $z=0$ in (\ref{120}), we immediately obtain $A_0(x,y)$ $=f(x,y,0)$. The proof of the assertion (\ref{shha}) of Theorem \ref{aprodpos} is now completed by  substituting (\ref{ZddA}) back into (\ref{120}).
\end{proof}
 \section{Summation formula for trivariate $q$-polynomials}
In this section, we give the summation formula for trivariate $q$-polynomials by direct computation and derive summation formulae for the second Hahn polynomials $\psi_n^{(a)}(x,y|q)$.  
\begin{thm}
\label{thm_1}
 Let the trivariate $q$-polynomials $F_n(x,y,z;q)$  be defined as in (\ref{3.6})
\begin{multline}
 F_{k+l} (x,\xi ,\zeta ;q)= 
 \sum_{n =0}^{k }\sum_{ r=0}^{l }\begin{bmatrix} k\\n \end{bmatrix}_{q}\begin{bmatrix} l\\r  \end{bmatrix}_{q} (-1)^{ n+r} q^{- \binom{n+1}{2}-r(n+l+1)-(k+l)(n+r)} 
P_{n+r} (\xi \ominus _{q}y,\zeta \ominus
_{q}z|q) \; F_{k+l-n-r} (x,y,z;q). \label{thm_1e}
\end{multline}
and
\begin{equation}
 F_{l} (x,\xi ,\zeta ;q)= 
\sum_{ r=0}^{l }\begin{bmatrix} l\\r  \end{bmatrix}_{q} (-1)^{ r} q^{-r(2l+1)} 
P_{r} (\xi \ominus _{q}y,\zeta \ominus
_{q}z|q) \; F_{l-r} (x,y,z;q). \label{thm_1ee}
\end{equation}
\end{thm}
Putting $\zeta =z$ in the last equation we establish the following corollary.
\begin{cor}
\begin{equation}
  F_{k+l} (x,\xi ,z ;q)=\sum_{n =0}^{k }\sum_{ r=0}^{l }\begin{bmatrix} k\\n \end{bmatrix}_{q}\begin{bmatrix} l\\r  \end{bmatrix}_{q} (-1)^{ n+r} q^{- \binom{n+1}{2}-r(n+l+1)-(k+l)(n+r)} 
  (\xi \ominus _{q}y)^{n+r} \; F_{k+l-n-r} (x,y,z;q)
.\label{corc}
\end{equation}
For $k=0$, we have:
\begin{equation}
  F_{l} (x,\xi ,z ;q)= \sum_{ r=0}^{l }\begin{bmatrix} l\\r  \end{bmatrix}_{q} (-1)^{  r} q^{ -r(2l+1)} 
  (\xi \ominus _{q}y)^{r} \; F_{l-r} (x,y,z;q)
.
\end{equation}
\end{cor}

\begin{proof}
 By replacing $t$ by $u\oplus _{q}t$ in  (\ref{3.6})
we get 
\begin{equation}
e_{q}\left[y\left( u\oplus _{q}t\right)\right]E_{q}\left[ -x\left( u\oplus
_{q}t\right) \right] {E}_{q} \left[-z\left( u\oplus _{q}t\right)\right]=\sum_{n=0}^{\infty }F_n(x,y,z;q)\frac{(-1)^nq^{({}^n_2)}\left( u\oplus _{q}t\right) ^{n}}{(q;q)_n}.  \label{3.13}
\end{equation}%
We now apply to the r.h.s of (\ref{3.13}) the identity

\begin{equation}
\sum_{j=0}^{+\infty }F(j)\frac{(x\oplus _{q}y)^{j}}{(q;q)_j}%
=\sum_{j,s=0}^{+\infty }F(j+s)\frac{q^{\binom{s}{2}}x^{j}y^{s}}{%
(q;q)_j(q;q)_s}.  \label{3.14}
\end{equation}%
satisfied by the JHC $q$-addition. So that (\ref{3.13}) becomes 
\begin{equation}
E_{q} \left[-x\left( u\oplus _{q}t\right)\right]=\frac{1}{%
e_{q}\left[y\left( u\oplus _{q}t\right)\right]E_{q }\left[ -z\left( u\oplus
_{q}t\right) \right] }\sum_{k,l=0}^{\infty }\frac{(-1)^{k+l} u^{k}t^{l}}{%
(q;q)_k(q;q)_l}\;q^{ \binom{k+l}{2}+\binom{l}{2}} F_{k+l} (x,y,z;q).  \label{3.15}
\end{equation}%
Note that the r.h.s of (\ref{3.15}) is independent of  variables $y$ and $z$ so that we
can write for any two variables $\xi ,\zeta $ the following equality 
\begin{equation}
\sum_{k,l=0}^{\infty }\frac{(-1)^{k+l} u^{k}t^{l}}{%
(q;q)_k(q;q)_l}\;q^{ \binom{k+l}{2}+\binom{l}{2}}\,F_{k+l} (x,\xi ,\zeta ;q)=\lambda
_{u,t}\left( \xi ,\zeta ;y,z;q\right) \sum_{k,l=0}^{\infty }\frac{(-1)^{k+l} u^{k}t^{l}}{%
(q;q)_k(q;q)_l}\;q^{ \binom{k+l}{2}+\binom{l}{2}} F_{k+l} (x,y,z;q).  \label{3.16}
\end{equation}%
where 
\begin{equation}
\lambda _{u,t}\left( \xi ,\zeta ;y,z;q\right) :=\frac{e_{q}\left[\xi \left(
u\oplus _{q}t\right) \right]E_{q }\left[ -\zeta \left( u\oplus _{q}t\right)\right] }{e_{q}\left[y\left( u\oplus _{q}t\right) \right]E_{q }\left[ -z\left(
u\oplus _{q}t\right)  \right] }.  \label{3.17}
\end{equation}%
By using the  rules in (\ref{2.12}), on can check that the quantity  (\ref{3.17}) also reads
\begin{equation}
\lambda _{u,t} \left( \xi ,\zeta ;y,z;q\right) =e_{q}\left[ \left( \xi \ominus
_{q}y\right) \left( u\oplus _{q}t\right) \right]E_{q }\left[ \left( \zeta
\ominus _{q}z\right) \left( u\oplus _{q}t\right)  \right].  \label{3.19}
\end{equation}%
On another hand, the r.h.s  of (\ref{3.19}) coincides with the generating function 
\begin{equation}
\sum_{n,r=0}^{\infty }\frac{u^{n}t^{r}}{(q;q)_n(q;q)_r}\;q^{\binom{r}{2}}%
P_{n+r} (\xi \ominus _{q}y,\zeta \ominus
_{q}z)  \label{3.20}
\end{equation}%
involving the Cauchy   polynomials.
Summarizing the above calculations in  (\ref{3.16})-(\ref{3.20}), we arrive at the sum
\begin{equation*}
\sum_{n,r=0}^{\infty }\frac{u^{n}t^{r}}{(q;q)_n(q;q)_r}\;q^{\binom{r}{2}}%
P_{n+r} (\xi \ominus _{q}y,\zeta \ominus
_{q}z|q)\sum_{k,l=0}^{\infty }\frac{(-1)^{k+l} u^{k}t^{l}}{%
(q;q)_k(q;q)_l}\;q^{ \binom{k+l}{2}+\binom{l}{2}} F_{k+l} (x,y,z;q)
\end{equation*}%
\begin{equation}
=\sum_{k,l=0}^{\infty }\frac{(-1)^{k+l} u^{k}t^{l}}{%
(q;q)_k(q;q)_l}\;q^{ \binom{k+l}{2}+\binom{l}{2}}\,F_{k+l} (x,\xi ,\zeta ;q).  \label{3.21}
\end{equation}%
Next, applying the series manipulation \cite[p.100]{SM}
\begin{equation}
\sum_{p=0}^{\infty }\sum_{s=0}^{\infty }A\left( p,s\right)
=\sum_{p=0}^{\infty }\sum_{s=0}^{p}A\left( s,p-s\right).  \label{3.22}
\end{equation}%
 to the l.h.s of (\ref{3.21}), we obtain that%
\begin{equation*}
\sum_{k,l=0}^{\infty }\sum_{n =0}^{k }\sum_{ r=0}^{l }\frac{(-1)^{k+l-n-r} u^{k}t^{l}q^{ \binom{k+l-n-r}{2}+\binom{l-r}{2}+\binom{r}{2}}}{(q;q)_{k-n}(q;q)_{l-r}(q;q)_n(q;q)_r} %
P_{n+r} (\xi \ominus _{q}y,\zeta \ominus
_{q}z|q) \; F_{k+l-n-r} (x,y,z;q)
\end{equation*}%
\begin{equation}
=\sum_{k,l=0}^{\infty }\frac{(-1)^{k+l} u^{k}t^{l}}{%
(q;q)_k(q;q)_l}\;q^{ \binom{k+l}{2}+\binom{l}{2}}\,F_{k+l} (x,\xi ,\zeta ;q).  \label{3.23}
\end{equation}%
By equating terms with $\displaystyle u^{k}t^{l}/(q;q)_k(q;q)_l$ and using the
simple combinatorial  fact 
we arrive at the following result%
\begin{equation*}
  F_{k+l} (x,\xi ,\zeta ;q)= \sum_{n =0}^{k }\sum_{ r=0}^{l }\begin{bmatrix} k\\n \end{bmatrix}_{q}\begin{bmatrix} l\\r  \end{bmatrix}_{q} (-1)^{ n+r} q^{- \binom{n+1}{2}-r(n+l+1)-(k+l)(n+r)} 
P_{n+r} (\xi \ominus _{q}y,\zeta \ominus
_{q}z|q) \; F_{k+l-n-r} (x,y,z;q).
\end{equation*}%
\end{proof}

As particular cases, by putting $y=\xi=ax,\,z=y$ and $\zeta=\xi$  in   Theorem \ref{thm_1}, we ave the  connection formulae for second Hahn polynomials $\psi_n^{(a)}(x,y|q)$  
as
\begin{cor}
\begin{eqnarray}
 \psi_{k+l}^{(a)} (x,\xi|q)= 
 \sum_{n =0}^{k }\sum_{ r=0}^{l }\begin{bmatrix} k\\n \end{bmatrix}_{q}\begin{bmatrix} l\\r  \end{bmatrix}_{q}  q^{ \binom{n+r}{2}- \binom{n+1}{2}-r(n+l+1)-(k+l)(n+r)} 
(\xi \ominus
_{q}y)^{n+r} \psi_{k+l-n-r}^{(a)} (x,y|q)\label{cor_31}
\end{eqnarray}
and
\begin{equation}
  \psi_{l}^{(a)} (x,\xi|q)= 
  \sum_{ r=0}^{l } \begin{bmatrix} l\\r  \end{bmatrix}_{q}  q^{ \binom{r}{2}- r(2l+1)} 
(\xi \ominus
_{q}y)^{r} \psi_{l-r}^{(a)} (x,y|q).\label{cor_32}
\end{equation}
\end{cor}
%\section{Summation formula for the  product of trivariate $q$-polynomials}
%In this section, we give and prove the summation formula for the  product of trivariate $q$-polynomials
 \begin{thm} 
 \label{thm_3}The following summation formula for the
product of  trivariate $q$-polynomials
\begin{align}
\label{Thm3}
   F_n(x,\xi,\zeta;q)F_r(X, \Omega, U;q) & = \sum_{ k,m=0}^{n,r} \begin{bmatrix}n\\k  \end{bmatrix}_{q} \begin{bmatrix} r\\m  \end{bmatrix}_{q}(-1)^{k+m}q^{\binom{k+1}{2}  + \binom{m+1}{2}  -mr-nk}   P_{k}(\xi\ominus_q  \zeta, y\ominus z)F_{n-k }(x,y,z;q) \cr
& \times\;  P_{m}(\Omega\ominus_q  U, Y\ominus_q  Z) F_{ r-m}(X, Y, Z;q).   
\end{align}
 holds true.
 \end{thm}
\begin{proof}  From the generating function (\ref{3.6}), we have%
\begin{equation}
e_{q}(y t)E_{q}(-xt){E}_{q}(-z t)e_{q}(Y T)E_{q}(-XT){E}_{q}(-Z T)=\sum_{n,r=0}^{\infty }\frac{(-1)^{n+r}q^{({}^n_2)+({}^r_2)} t^{n}T^r}{(q;q)_n(q;q)_r}F_n(x,y,z;q)F_r(X, Y, Z;q).   \label{3.27}
\end{equation}%
Replacing in  (\ref{3.27}) $y$ by $\xi ,$ $z$ by $\zeta ,$ $Y$ by $%
\Omega $ and $Z$ by $U,$ we get
\begin{equation}
e_{q}(\xi t)E_{q}(-xt){E}_{q}(-\zeta t)e_{q}(\Omega T)E_{q}(-XT){E}_{q}(-U T)=\sum_{n,r=0}^{\infty }\frac{(-1)^{n+r}q^{({}^n_2)+({}^r_2)} t^{n}T^r}{(q;q)_n(q;q)_r}F_n(x,\xi,\zeta;q)F_r(X, \Omega, U;q).   \label{3.28}
\end{equation}
By replacing 
\begin{eqnarray*}
  E_{q}(-xt) E_{q}(-XT) =\frac{1}{e_{q}(y t) {E}_{q}(-z t)e_{q}(Y T){E}_{q}(-Z T)}\sum_{n,r=0}^{\infty }\frac{(-1)^{n+r}q^{({}^n_2)+({}^r_2)} t^{n}T^r}{(q;q)_n(q;q)_r}F_n(x,y,z;q)F_r(X, Y, Z;q). 
\end{eqnarray*}
in  the l.h.s. of  (\ref{3.28}) and using (\ref{2.12})
%%% 
 one gets, after expanding the exponentials in series, the following
\begin{align}
 & \sum_{n,r=0}^{\infty }\frac{(-1)^{n+r}q^{({}^n_2)+({}^r_2)} t^{n}T^r}{(q;q)_n(q;q)_r}F_n(x,\xi,\zeta;q)F_r(X, \Omega, U;q) \cr
&\qquad\qquad =\frac{e_{q}(\xi t) {E}_{q}(-\zeta t)e_{q}(\Omega T) {E}_{q}(-U T)}{e_{q}(y t) {E}_{q}(-z t)e_{q}(Y T) {E}_{q}(-Z T)}\sum_{n,r=0}^{\infty }\frac{(-1)^{n+r}q^{({}^n_2)+({}^r_2)} t^{n}T^r}{(q;q)_n(q;q)_r}F_n(x,y,z;q)F_r(X, Y, Z;q)\cr
& \qquad\qquad
=\frac{e_q\Big( (\xi  \ominus_q  \zeta)t\Big)}{ e_q\Big( (y  \ominus_q  z)t\Big)}\frac{ e_q\Big( (\Omega  \ominus_q  U)T\Big)}{e_q\Big( (Y   \ominus_q  Z)T\Big)}  
\sum_{n,r=0}^{\infty }\frac{(-1)^{n+r}q^{({}^n_2)+({}^r_2)} t^{n}T^r}{(q;q)_n(q;q)_r}F_n(x,y,z;q) F_r(X, Y, Z;q)
\cr
& \qquad\qquad=\sum_{n,k=0}^{\infty}  \frac{(-1)^{n }q^{({}^n_2)} t^{n+k} }{(q;q)_n(q;q)_k}P_{k}(\xi\ominus_q  \zeta, y\ominus z)F_n(x,y,z;q) \cr
&\qquad\qquad \times\;\sum_{r,m=0}^{\infty}  \,   \frac{(-1)^{ r}q^{ ({}^r_2)} T^{r+m} }{(q;q)_r(q;q)_m}P_{m}(\Omega\ominus_q  U, Y\ominus_q  Z) F_r(X, Y, Z;q).   
\end{align}
Finally, by replacing $n$ by $n-k$ and $r$ by $r-m$ in the r.h.s. of the last relation, the proof is completed.  
%\begin{align*}
% & 
%\sum_{n,r=0}^{\infty }\frac{(-1)^{n+r}q^{({}^n_2)+({}^r_2)} t^{n}T^r}{(q;q)_n(q;q)_r}F_n(x,\xi,\zeta;q)F_r(X, \Omega, U;q) \cr
%&\qquad\qquad=\sum_{n =0}^{\infty}\sum_{ k=0}^n  \frac{(-1)^{n-k }q^{({}^{n-k }_2)} t^{n } }{(q;q)_{n-k }(q;q)_k}P_{k}(\xi\ominus_q  \zeta, y\ominus z)F_{n-k }(x,y,z;q) \cr
%&\qquad\qquad\times\;\sum_{r =0}^{\infty}  \sum_{m=0}^r  \frac{(-1)^{ r-m}q^{ ({}^{ r-m}_2)} T^{r} }{(q;q)_{ r-m}(q;q)_m}P_{m}(\Omega\ominus_q  U, Y\ominus_q  Z) F_{ r-m}(X, Y, Z;q).   (3.29)
%\end{align*}
%Finally
%\begin{align}
%   F_n(x,\xi,\zeta;q)F_r(X, \Omega, U;q) & = \sum_{ k,m=0}^{n,r} \begin{bmatrix}n\\k  \end{bmatrix}_{q} \begin{bmatrix} r\\m  \end{bmatrix}_{q}(-1)^{k+m}q^{({}^{k+1 }_2)+({}^{ m+1}_2)-mr-nk}   P_{k}(\xi\ominus_q  \zeta, y\ominus z)F_{n-k }(x,y,z;q) \cr
%& \times\;  P_{m}(\Omega\ominus_q  U, Y\ominus_q  Z) F_{ r-m}(X, Y, Z;q).   (3.29)
%\end{align}
\end{proof}
\section{Particular cases}
\begin{enumerate}

\item Putting $y=\xi=ax,\,z=y$ and $\zeta=\xi$  in   Theorem \ref{thm_1}, leads to a connection formulae for second Hahn polynomials $\psi_n^{(a)}(x,y|q)$  
\begin{eqnarray}
 \psi_{k+l}^{(a)} (x,\xi|q)= 
 \sum_{n =0}^{k }\sum_{ r=0}^{l }\begin{bmatrix} k\\n \end{bmatrix}_{q}\begin{bmatrix} l\\r  \end{bmatrix}_{q}  q^{ \binom{n+r}{2}- \binom{n+1}{2}-r(n+l+1)-(k+l)(n+r)} 
P_{n+r}(\xi,y) \psi_{k+l-n-r}^{(a)} (x,y|q).
\end{eqnarray}
and
\begin{equation}
  \psi_{l}^{(a)} (x,\xi|q)= 
  \sum_{ r=0}^{l } \begin{bmatrix} l\\r  \end{bmatrix}_{q}  q^{ \binom{r}{2}- r(2l+1)} 
P_r(\xi, y) \psi_{l-r}^{(a)} (x,y|q).
\end{equation}
\item   Putting $y=\xi=ax,\, Y=\Omega=aX,\,z=y,\,\zeta=\xi$ and $U=\Omega$ in Theorem \ref{thm_3}, we get
\begin{align}
\label{Thm3}
   \psi_n^{(a)}(x,\xi|q)\psi_r^{(a)}(X, \Omega|q) & = \sum_{ k,m=0}^{n,r} \begin{bmatrix}n\\k  \end{bmatrix}_{q} \begin{bmatrix} r\\m  \end{bmatrix}_{q}(-1)^{k+m}q^{\binom{k+1}{2} + \binom{m+1}{2}-mr-nk}   P_{k}(ax\ominus_q  \xi, ax\ominus y)\;\psi_{n-k }^{(a)}(x,y|q) \cr
& \times\;  P_{m}(aX\ominus_q  \Omega, aX\ominus_q  Y)\, \psi_{ r-m}^{(a)}(X, Y|q).   
\end{align}

\end{enumerate}

\end{CJK*}
\end{document}